\documentclass[12pt]{article}
\pdfoutput=1
\usepackage{euscript,amsmath,amsthm,amsfonts,amssymb}
\usepackage{mathtools}
\usepackage{subfig}
\usepackage{float}
\usepackage[margin=1in]{geometry}
\usepackage{graphicx}
\usepackage{outlines}
\usepackage[dvipsnames]{xcolor}
\usepackage[noautoscale]{youngtab}
\usepackage{relsize}
\usepackage{jheppub}

\newcommand{\RNum}[1]{\uppercase\expandafter{\romannumeral #1\relax}}

\newcommand{\be}{\begin{equation}}
\newcommand{\ee}{\end{equation}}
\newcommand{\bea}{\begin{eqnarray}}
\newcommand{\eea}{\end{eqnarray}}

\title{\boldmath Level-rank duality of knot and link invariants}
\author{Howard J. Schnitzer}
\affiliation{Martin Fisher School of Physics, Brandeis University, Waltham, Massachusetts 02453, USA}
\preprint{BRX-TH-6684}
\emailAdd{schnitzr@brandeis.edu}
\abstract{A number of results for the level-rank duality of $G(N)_{K}\longleftrightarrow G(K)_N$ Chern-Simons theory are summarized, with emphasis on the applications to knot and link invariants. Explicit examples for $SU(2)_K\longleftrightarrow SU(K)_2$ illustrate general results.

A criterion to distinguish torus knots and links from hyperbolic knots and links, based on tables constructed by Kaul for one and two strand invariants, is presented. Symmetries of hyperbolic knot and link invariants are discussed. 

The level-rank duality of torus knot and link invariants of minimal models is examined.

}
\begin{document}
\Yboxdim6pt
\maketitle
\flushbottom

\section{Introduction}
Group-level duality of $G(N)_K$ Chern-Simons theory \cite{r1,r2,r3,r4,r5} has provided a number of insights into the structure of topological theories. Properties of pairs of theories are tightly constrained by duality maps $G(N)_{K}\longleftrightarrow G(K)_N$, where $N$ is the dimension of the defining representation of $SU(N)$ and $K$ the dimension of the associated affine Lie algebra, with analogous definitions for $SO(N)_K,Sp(N)_{K}$, and the exceptional groups.

There are dualities of primary fields, conformal dimensions, braid matrices, skein relations, and knot and link invariants \cite{r3,r4,r5}. Specifically every knot in $S^3$ is either a torus knot, a satellite knot, or a hyperbolic knot \cite{r6}. A torus knot is a knot that can be embedded without crossings on the surface of a torus in $S^3$. A satellite knot is a knot that can be embedded in the regular neighborhood of another knot in $S^3$ \cite{r7,r8}. Every non-split alternating link that is not a torus link is hyperbolic. An example of a hyperbolic knot is the figure 8 knot, and that of hyperbolic links are the Whitehead link and the Borromean link.

Chern-Simons knot theory was pioneered by Witten \cite{r10,r11}. The level-rank dualities of knot and link invariants are presented in section \ref{sec:2}; with explicit examples for $SU(2)_K \longleftrightarrow SU(K)_2$ given in subsection \ref{sub23}. Section \ref{sec:3} considers hyperbolic knots and links. A conjectured criterion for distinguishing torus knots and links from hyperbolic knots and links is based on systematic features of the tables of Kaul \cite{r13}. This criterion is satisfied by the Whitehead link and the Borromean link. In section \ref{sec:4}, the level-rank duality of minimal models is discussed.

In appendix \ref{sec:A} we give evidence for symmetries of $U(N)$ hyperbolic knot and link invariants.

\section{$G(N)_{K}\longleftrightarrow G(K)_N$}\label{sec:2}
The conformal dimension of a $\lambda$-primary of $G(N)_K$ \cite{r3},
\be
h_\lambda =\frac{1}{2}Q_{\lambda}(N) / (K+g)
\ee
is related to
\be
\tilde{h}_{\tilde{\lambda}} =\frac{1}{2}Q_{\tilde{\lambda}}(K)/ (N+\tilde{g})
\ee
where $g$ is the dual Coxeter number for $G(N)$, and $\tilde{g}$ is the same for $G(K)$, with $Q_{\lambda}(N)$ the quadratic Casimir operator for representation $\lambda$, and similarly $Q_{\tilde{\lambda}}(K)$ for representation $\tilde{\lambda}$. The Young tableau for representation $\lambda$ is related to $\tilde{\lambda}$, which is the transpose of the tableau $\lambda$.

The conformal dimensions of $SU(N)_K$ satisfy the duality relation \cite{r4,r5}
\be
h(a)+\tilde{h}(\tilde{a})=\frac{1}{2}r(a)\left[1-\frac{r(a)}{2NK}\right]
\ee
where $r(a)$ is the number of boxes in the Young tableau for representation $a$. Also relevant is the root of unity
\be
q=e^{\frac{2 \pi i}{N+K}}
\ee
for $SU(N)_K$.

The modular transformation and fusion matrices of $SU(N)_K$ and $SU(K)_N$ satisfy \cite{r4}
\be\label{eq:25}
S_{ab}=\sqrt{\frac{K}{N}}\exp\left[\frac{-2\pi i r(a)r(b)}{NK}\right]\tilde{S}^*_{\tilde{a}\tilde{b}}
\ee
and
\be
N_{ab}^c=\tilde{N}_{\tilde{a}\tilde{b}}^{\sigma^\Delta(\tilde{c})}
\ee
for
\be\label{eq:27}
\Delta=\frac{r(a)+r(b)-r(c)}{N}
\ee
In general the nuber of boxes of $c$ may be less than the sum of those of $a$ and $b$, because $c$ has been ``reduced". In that case $c$ is dual to the representation $\sigma^\Delta(\tilde{c})$ which is ``cominamally equivalent'' to $\tilde{c}$. For example, for $SU(2)_K$ the fusion matrix is
\be
N_{ab}^c, \quad \text{where } c=a+b\mod K
\ee
If $(a+b)<K$, then $c=(a+b)$. But if $(a+b)\geq K$, then $c=(a+b)-K\leq (a+b)$.

In general, to construct the Young tableau for $\sigma^\Delta(\tilde{c})$, there is a specific procedure.
\begin{enumerate}
    \item First transpose the reduced diagram c. The row lengths (not necessarily reduced) transform to $\tilde{l}_i(\tilde{c})$.
    \item Reduce this diagram.
    \item Add a row of length $N$ boxes to the top of this diagram.
    \item Reduce this diagram.
    \item Repeat these steps $(\Delta-1)$ times.
\end{enumerate}
The final reduced diagram will have
\be
r(c)+N\Delta-K\tilde{l}_{K-\Delta}
\ee
boxes. That is
\be
(a,b,c) \longrightarrow (\tilde{a},\tilde{b},\sigma^\Delta(\tilde{c}))
\ee
for the transform of the fusion matrix.

\subsection{Knot dualities}\label{sub21}
With these ingrediants one can compute $S_{0a}/S_{00}$ and it's level-rank dual to show that \cite{r3}
\be
\langle\text{unknot};a\rangle_{G(N)_K}=\langle\text{unknot};\tilde{a}\rangle_{G(K)_N}
\ee
More generally for knots that can be untied by skein relations (not effective in all cases)\cite{r3}
\be\label{eq:212}
\langle \mathcal{K},\lambda\rangle_{G(N)_{K}}=\langle \tilde{\mathcal{K}},\tilde{\lambda}\rangle_{G(K)_N}
\ee
in vertical framing, where $\tilde{\mathcal{K}}$ is the mirror image of knot $\mathcal{K}$.

\subsection{Link dualities}\label{sub22}
The duality \cite{r5} for two strands is
\be\label{eq:213}
\mathcal{L}(a,b;\{n_i\})_{G(N)_K} = e^{-i\pi \Phi(a,b)\sum_in_i}\tilde{\mathcal{L}}(\tilde{a},\tilde{b};\{-n_i\})_{G(K)_N}
\ee
for vertical framing, where $\sum n_i= \text{even}$, and each $n_i$ is the sum of crossings of the neighboring braids, and 
\be\label{eq:214}
\Phi(a,b)=\frac{r(a)r(b)}{NK}
\ee
for $SU(N)_K$. More generally, in vertical framing
\be\label{eq:215}
\mathcal{L}(\{a_i\})_{G(N)_K}=e^{i\pi \sum_iw(i,i)r(a_i)}e^{-i\pi\sum_{i,j}w(i,j)\Phi(a_i,a_j)}\tilde{\mathcal{L}}(\{\tilde{a}_i\})_{G(K)_N}
\ee
The mirror image of link $\mathcal{L}$ is $\tilde{\mathcal{L}}$, and $w(i,j)$ is the sum of crossings between components $i$ and $j$.

It is to be emphasized that the derivation of \eqref{eq:212} and \eqref{eq:215} require the reduction to planar graphs and planar pentagons. Is this a coincidence, since \emph{hyperbolic} knots and links involve non-planar polyhedra? Or is there a deeper connection \cite{r7,r8}?

\subsection{$SU(2)_K\longleftrightarrow SU(K)_2$}\label{sub23}
In this subsection we present some examples of link invariants to illustrate the dualities discusses in subsections \ref{sub21} and \ref{sub22}.

\paragraph{Hopf link: $2^2_1$ in Rolfson notation}

The Hopf link for two strands belonging to representation $j_1$ and $j_2$ in $SU(2)_K$ is \cite{r11}
\be
|j_1,j_2\rangle=S_{j_1,j_2}/S_{00},
\ee
where $j_1=\frac{a_1}{2}$ and $j_2=\frac{a_2}{2}$, with $a_1$ and $a_2$ the number of boxes in the representation. From equations \ref{eq:25} and \eqref{eq:214}
\be
\begin{split}
|j_1,j_2\rangle &= \exp \left[-\frac{2\pi ir(a_1)r(a_2)}{2K}\right]\frac{\tilde{S}^*_{\tilde{a}_1\tilde{a}_2}}{\tilde{S}^*_{00}}\\
&=\exp \left[-2\pi i\Phi(a_1,a_2)\right]\frac{\tilde{S}^*_{\tilde{a}_1\tilde{a}_2}}{\tilde{S}^*_{00}}\\
&=\exp \left[-2\pi i\Phi(a_1,a_2)\right]|\tilde{j}_1,\tilde{j}_2\rangle
\end{split}
\ee
in agreement with \eqref{eq:213}, with $
\sum_in_i=2$.

\paragraph{Sum of two Hopf links}
\cite{r11,r12}
\be
\begin{split}
|2^2_1+2^2_1\rangle &= |j_1,j_2,j_3\rangle\\
&=\frac{S_{j_1j_2}S_{j_3j_2}}{S_{0j_2}S_{00}}
\end{split}
\ee
using equations \eqref{eq:25} and \eqref{eq:214}, we find
\be
|j_1j_2j_3\rangle =\exp\left(-\frac{2\pi i}{2K}r(a_2)\left[r(a_1)+r(a_2)\right]\right)|\tilde{j}_1\tilde{j}_2\tilde{j}_3\rangle
\ee
in agreement with \eqref{eq:25} and \eqref{eq:214}, with $w(1,2)=w(2,3)=2$.

\paragraph{Link $6^3_3$}

For fixed framing \cite{r12},
\be
|j_1j_2j_3\rangle =\sum_m\exp\left[h_m-h_{j_1}-h_{j_3}\right]\frac{N_{a_1a_3}^mS_{mj_2}}{S_{00}^3}
\ee
From \eqref{eq:25}-\eqref{eq:27} and \eqref{eq:214}, the result is
\be
|j_1j_2j_3\rangle =\exp\left[-i\pi\left(\Phi(a_1,a_2)+\Phi(a_2,a_3)\right)\right]\exp\left[-2\pi i\left(\tilde{h}_{a_1+a_3}-\tilde{h}_{a_1}-\tilde{h}_{a_3}\right)\right]|\tilde{j}_1\tilde{j}_2\tilde{j}_3\rangle.
\ee
for the special case where $a_1+a_3<K$, and more generally
\be
|j_1j_2j_3\rangle =\exp\left[-i\pi\left(\Phi(a_1,a_2)+\Phi(a_2,a_3)\right)\right]\exp\left[-2\pi i\left(\tilde{h}_{\sigma^\Delta}(\tilde{c})-\tilde{h}_{a_1}-\tilde{h}_{a_3}\right)\right]|\tilde{j}_1\tilde{j}_2\tilde{j}_3\rangle.
\ee
This satisfies \eqref{eq:213} up to framing, where \eqref{eq:213} is for vertical framing. Vertical framing can be restored using
\be
T|m\rangle = e^{2\pi ih_m}|m\rangle
\ee
where $T$ is a modular transformation operator.

These examples of torus links illustrate that they transform under level-rank duality in the expected way.

\section{Hyperbolic knots and links: $SU(2)_K$}\label{sec:3}
The purpose of this chapter is to present evidence for a criterion which distinguishes hyperbolic knots and links from torus knots and links. Begin with the hyperbolic Whitehead link with two strands and the Borromean link with three strands.
\paragraph{Whitehead link: $5_1^2$}
The result of the construction of Kaul \cite{r12,r13} is
\be\label{eq:41}
\begin{split}
  C_{5_1^2}(j_1j_2)&=\left[2j_1+1\right]^2\left[2j_2+1\right]\sum_{m,n,p}\lambda^{-1}_{p_1,-}(j_1,j_2)\lambda_{p_2,+}(j_1,j_2) \lambda^{-1}_{n_1,+}(j_1,j_2) \lambda^{-1}_{m_1,-}(j_1,j_2) \lambda^{-1}_{m_2,+}(j_1,j_2) \\
  &\times a_{(0,p)}\begin{pmatrix}
  j_1 & j_1\\ 
  j_2 & j_2\\
  j_1 & j_1
\end{pmatrix}
a_{(n,p)}\begin{pmatrix}
  j_1 & j_2\\ 
  j_1 & j_1\\
  j_2 & j_1
\end{pmatrix}
a_{(n,m)}\begin{pmatrix}
  j_1 & j_2\\ 
  j_1 & j_1\\
  j_2 & j_1
\end{pmatrix}
a_{(0,m)}\begin{pmatrix}
  j_1 & j_1\\ 
  j_2 & j_2\\
  j_1 & j_1
\end{pmatrix}.
\end{split}
\ee
where $a_{(n,p)}$'s are duality transformations acting on 6-point conformal blocks on $S^2$, and the $\lambda$'s are phases which these blocks pick up under the action of braid generators. Equivalently, the $a_{(n,p)}$ are related to the quantum $6j$ symbols. Explicit calculations verify that there is no $SU(2)_K \longleftrightarrow SU(K)_2$ duality for \eqref{eq:41}. The bracket $[2j+1]$ is the quantum dimension
\be
[x]=\frac{q^{\frac{x}{2}}-q^{-\frac{x}{2}}}{q^{\frac{1}{2}}-q^{-\frac{1}{2}}}
\ee
where $q=\exp[2\pi i/(K+N)]$ for $SU(N)_K$. See also \cite{r20}.

\paragraph{Borromean link: $6_2^3$}
The construction of Kaul \cite{r12,r13} gives
\be\label{eq:43}
\begin{split}
  C_{6_2^3}(j_1j_2)&=\left[2j_1+1\right]\left[2j_2+1\right]\left[2j_3+1\right]\sum_{l,m,n,p,q}\lambda_{l_1,-}(j_1,j_2)\lambda_{l_2,-}(j_1,j_3) \lambda^{-1}_{m_1,-}(j_2,j_3)\lambda^{-1}_{n_2,+}(j_1,j_2)\\
  &\times\lambda_{p_0,+}(j_1,j_2)\lambda^{-1}_{p_1,-}(j_1,j_3)\lambda^{-1}_{q_1,}(j_1,j_2)\lambda_{q_2,-}(j_2,j_3) \\
  &\times a_{(0,l)}\begin{pmatrix}
  j_2 & j_2\\ 
  j_1 & j_1\\
  j_3 & j_3
\end{pmatrix}
a_{(m,l)}\begin{pmatrix}
  j_2 & j_1\\ 
  j_2 & j_3\\
  j_1 & j_3
\end{pmatrix}
a_{(m,n)}\begin{pmatrix}
  j_2 & j_1\\ 
  j_3 & j_2\\
  j_1 & j_3
\end{pmatrix}\\
&\times a_{(p,n)}\begin{pmatrix}
  j_2 & j_1\\ 
  j_3 & j_1\\
  j_2 & j_3
\end{pmatrix}
a_{(p,q)}\begin{pmatrix}
  j_1 & j_2\\ 
  j_1 & j_3\\
  j_2 & j_3
\end{pmatrix}
a_{(0,q)}\begin{pmatrix}
  j_1 & j_1\\ 
  j_2 & j_2\\
  j_3 & j_3
\end{pmatrix}.
\end{split}
\ee
Again there is no level-rank duality for \eqref{eq:43}.

\paragraph{Knot invariants}
Table \RNum{2}A of Kaul \cite{r13} provides a list of knot invariants where all knots carry spin $j$ of $SU(2)_K$. The torus knots listed are $3_1$, $5_1$, and $7_1$, which do not have factors of $a_{ml}$. All other knots are hyperbolic and have one or more factors of $a_{ml}$.

\paragraph{Link invariants}
Table \RNum{2}B of Kaul \cite{r13} lists the invariants of two component links up to seven crossings for spins $j_1$ and $j_2$ on the component strands. The torus links $0_1$, $2_1$, $4_1$, and $6_1$ do not have factors of $a_{(m)(l)}$.

The consistent overview is that torus knots and links are distinguished from hyperbolic knots and links by the absence or presence of factors of $a_{ml}$ or $a_{(m)(l)}$. It would be interesting to prove this observation. See also \cite{r14,r15,r16,r17}. Brunnian links are also a class of hyperbolic links \cite{r35}.

\paragraph{Symmetry}
The torus knot and link invariants exhibit a level-rank duality, as discussed in sections \ref{sec:2}, but the hyperbolic invariants do not have this duality. Since the hyperbolic invariants are related to $SL(2,\mathbb{C})$ Chern-Simons theory, one may ask if there is an analogous symmetry for $SL(2,\mathbb{C})$ (see Appendix \ref{sec:A}). For other applications of $SU(N)_K$ Chern-Simons theory see \cite{r17,r18,r19,r24,r25,r26}.

\section{Minimal models}\label{sec:4}
Minimal models can be described by the coset \cite{r31,r32,r36}
\be\label{eq:441}
SU(2)_K \times SU(2)_1 /  SU(2)_{K+1}
\ee
Since the central charge of $SU(N)_K$ is
\be
c^{N,K}=\frac{(N^2-1)K}{N+K},
\ee
it is easy to show that
\be
c_{min}=1-\frac{6}{(K+2)(K+3)}
\ee
For low values of $K$, one has
\be
\begin{split}
    &K=1;\quad c=\frac{1}{2} \quad \text{Ising model,}\\
     &K=2;\quad c=\frac{7}{10} \quad \text{Tri-critical Ising Model,}\\
     & K=3;\quad c=\frac{4}{5} \quad \text{Potts model.}
\end{split}
\ee

The level-rank dual of \eqref{eq:441} is
\be\label{eq:445}
SU(K+1)_2 / \left[SU(K)_2 \times U(1)_2\right]
\ee
with central charge 
\be
c_{dual}=c_{min}
\ee
In this sense, the minimal models are level-rank self-dual.

Recall from \eqref{eq:215}, for link $L$, with representations $r_{1},r_2,\cdots r_n$ of the component knots, with $a(r_{i})$ denoting the number of boxes of representation $r_i$, so that
\be\label{eq:447}
\mathcal{L}(\{r_i\})_{SU(2)_K}=e^{i\pi \sum_iw(i,i)r(a(r_i)}e^{-i\pi\sum_{i,j}w(i,j)\Phi(r_i,r_j)}\tilde{\mathcal{L}}(\{\tilde{r}_i\})_{SU(K)_2}
\ee
The mirror image of Theorem 1 equation (25) of reference \cite{r30a} for minimal models is
\be
V_{(\bar{r}_1,\bar{s}_1)\cdots (\bar{r}_n,\bar{s}_n)}[L]=V^{(K)}_{\bar{r}_1\cdots \bar{r}_n}[\bar{L}]V^{(K+1)}_{\bar{s}_1\cdots \bar{s}_n}[L]V^{(1)}_{\bar{\epsilon}_1\cdots \bar{\epsilon}_n}[\bar{L}]
\ee
where $V^{(K)}_{\{\bar{r}_i\}}[\bar{L}]$, $V^{(K+1)}_{\{\bar{s}_i\}}[L]$, and $V^{(1)}_{\{\bar{\epsilon}_i\}}[\bar{L}]$, are the mirror image links for $SU(K)_2$,$SU(K+1)_2$, and $U(1)_2$ respectively, with $\bar{L}$ the mirror image link.

Analogous to \eqref{eq:447} one has
\be\label{eq:449}
\tilde{\mathcal{L}}(\{s_i\})_{SU(2)_{K+1}}=e^{-i\pi \sum_iw(i,i)r(a(s_i)}e^{i\pi\sum_{i,j}w(i,j)\Phi(s_i,s_j)}\mathcal{L}(\{\tilde{s}_i\})_{SU(K+1)_2}.
\ee
Further $V^{(1)}_{\epsilon_1,\cdots , \epsilon_n}(L)$ is the torus link invariant for $SU(2)_1$, where $\epsilon_i=1$ or 2 for $(r_i-s_i)$ even or equivalently of $U(1)_2$, and $\epsilon_i=2$ for the representation with one box or odd, with $\epsilon_i=1$ for the identity representation of $SU(2)_1$.

Putting this all together,
\be\label{eq:4410}
V^{(K)}_{r_1\cdots r_n}[L]V^{(K+1)}_{s_1\cdots s_n}[\bar{L}]V^{(1)}_{\epsilon_1\cdots \epsilon_n}[L]=(phases)\tilde{V}^{(K)}_{\tilde{r}_1\cdots \tilde{r}_n}[\bar{L}]\tilde{V}^{(K+1)}_{\tilde{s}_1\cdots \tilde{s}_n}[L]\tilde{V}^{(1)}_{\epsilon_1\cdots \epsilon_n}[\bar{L}]
\ee
where $\Phi$ is as in \eqref{eq:214}, and $w(i,j)$ is the sum of the crossings between components $i$ and $j$. Therefore $(phases)$ is obtained from the product of phases from \eqref{eq:447},\eqref{eq:449}, and for $SU(2)_1\longrightarrow U(1)_2$. Thus \eqref{eq:4410} relates the link invariants of \eqref{eq:441} to the mirror image link invariants of \eqref{eq:445}; explicitly expressing the self-duality of the minimal models.

An interesting question is how these invariants could emerge from lattice theories of minimal models? This subject is not as well developed as that obtained directly from rational conformal field theory described by the coset models. However special cases may be relevant. For example, non-abelian anyons on a spin $\frac{1}{2}$ honeycomb lattice has been studied by Kitaev \cite{r34}\footnote{We thank Djordje Radicevic for informing us of this reference}. The fusion and braid rules for the non-abelian anyons are described in Secs. 6 and 8 of that paper, where the boundary CFT is the Ising model. Once one has the braid group, one can construct torus knot and link invariants following Kaul and collaborators \cite{r14,r15,r16,r30a}. It is likely that to construct hyperbolic knots with links in lattice theories, one will require the lattice version of $6-j$ symbols as suggested by Sec. 3 above.

\section{Concluding remarks}
In this paper we discussed the level-rank duality for $G(N)_K\longleftrightarrow G(K)_N$ Chern-Simons theory, with examples which distinguish torus knot and link invariants from hyperbolic knot and link invariants. An open question is whether there is a similar duality for $SL(N,\mathbb{C})$ Chern-Simons theory.

\acknowledgments
We thank Jonathan Harper and Isaac Cohen-Abbo for their help in preparing the manuscript, and Djordje Radicevic for insightful comments.

\appendix

\section{A symmetry of knot and link invariants}\label{sec:A}
Ramadevi \cite{r22} has computed $U(N)_K$ polynomials for some non-torus knots, and two-component links, which generalize the $SU(2)_K$ results of Kaul \cite{r13}. The above invariants are polynomials in two-variables
\be\label{eq:A1}
q=\exp\left(\frac{2\pi i}{K+N}\right)
\ee
and 
\be\lambda=q^N.
\ee 
The knot polynomials have the symmetry relations for $U(N)$, $\lambda$ fixed \cite{r22}
\be\label{eq:A3}
V_{\yng(1)}[\mathcal{K}](q^{-1})=-V_{\yng(1)}[\mathcal{K}](q)
\ee
\be
V_{\yng(2)}[\mathcal{K}](q^{-1})=V_{\yng(1,1)}[\mathcal{K}](q).
\ee
Similarly the invariants for link polynomials satisfy \cite{r22}
\be
V_{(\yng(1),\:\yng(1))}[\mathcal{L}](q^{-1})=V_{(\yng(1),\:\yng(1))}[\mathcal{L}](q)
\ee
\be
V_{(\yng(1),\:\yng(2))}[\mathcal{L}](q^{-1})=-V_{(\yng(1),\:\yng(1,1))}[\mathcal{L}](q)
\ee
\be\label{eq:A7}
V_{(\yng(2),\:\yng(2))}[\mathcal{L}](q^{-1})=V_{(\yng(1,1),\:\yng(1,1))}[\mathcal{L}](q)
\ee
where $\yng(1),\yng(2),\yng(1,1)$ signify the defining, symmetric, and anti-symmetric representations respectively. There are no obvious symmetries associated with $\lambda$.

We propose the extension of \eqref{eq:A3}-\eqref{eq:A7} for knot and link invariants,
\be\label{eq:A8}
V_{a_1\cdots a_n}[\mathcal{L}](q^{-1})=(-1)^lV_{\tilde{a}_1\cdots \tilde{a}_n}[\mathcal{L}](q)
\ee
where $l$ is the sum of the number of boxes in representations $(a_1\cdots a_n)$ and $(\tilde{a}_1\cdots \tilde{a}_n)$ are the ``transposed'' representations.

Equations \eqref{eq:A3}-\eqref{eq:A8} have aspects of level rank symmetry, but $q\rightarrow q^{-1}$ is not the result of $(K\longleftrightarrow N)$ in \eqref{eq:A1}. We expect that \eqref{eq:A3}-\eqref{eq:A8} could be the result of symmetry operations involving representations of a non-compact group, such as $SL(N,\mathbb{C})$. Level-rank duality of the $U(N)_{K}$ WZW Chern-Simons theory was explored in \cite{r21}, which would apply to torus knots and links, but not \eqref{eq:A3}-\eqref{eq:A8}. See also \cite{r23,r24,r25,r26}.

\paragraph{Note added:}
 Equation \eqref{eq:A8} has been proven in equations (1.13) - (1.15) of \cite{r27} and (1.5) - (1.7 ) of \cite{r28}. We thank Professor Zhu for providing this information. We also thank Sergei Gukov for helpful comments and calling our attention to references \cite{r29}, \cite{r30}. A recent application of level-rank duality is \cite{r33}.

\bibliographystyle{JHEP}
\bibliography{main}

\providecommand{\href}[2]{#2}\begingroup\raggedright\begin{thebibliography}{10}

\bibitem{r1}
S.G.~Naculich and H.J.~Schnitzer, \emph{Duality relations between {SU(N)K} and
  {SU(K)N WZW} models and their braid matrices},
  \href{https://doi.org/https://doi.org/10.1016/0370-2693(90)90061-A}{\emph{Physics
  Letters B} {\bfseries 244} (1990) 235}.

\bibitem{r2}
A.~Kuniba and T.~Nakanishi, \emph{in \uppercase{P}roc. \uppercase{I}nt.
  \uppercase{C}oll. on modern quantum field theory, \uppercase{B}ombay},  1990.

\bibitem{r3}
S.~Naculich, H.~Riggs and H.~Schnitzer, \emph{Group-level duality in {WZW}
  models and chern-simons theory},
  \href{https://doi.org/https://doi.org/10.1016/0370-2693(90)90623-E}{\emph{Physics
  Letters B} {\bfseries 246} (1990) 417}.

\bibitem{r4}
E.J.~Mlawer, S.G.~Naculich, H.A.~Riggs and H.J.~Schnitzer, \emph{Group-level
  duality of {WZW} fusion coefficients and {Chern-Simons} link observables},
  \href{https://doi.org/https://doi.org/10.1016/0550-3213(91)90110-J}{\emph{Nuclear
  Physics B} {\bfseries 352} (1991) 863}.

\bibitem{r5}
S.G.~Naculich, H.A.~Riggs and H.J.~Schnitzer, \emph{Simple-current symmetries,
  rank-level duality, and linear skein relations for {Chern-Simons} graphs},
  \href{https://doi.org/https://doi.org/10.1016/0550-3213(93)90022-H}{\emph{Nuclear
  Physics B} {\bfseries 394} (1993) 445}.

\bibitem{r6}
W.P.~Thurston, \emph{Three dimensional manifolds, {Kleinian} groups and
  hyperbolic geometry},
  \href{https://doi.org/10.1090/s0273-0979-1982-15003-0}{\emph{Bulletin of the
  American Mathematical Society} {\bfseries 6} (1982) 357}.

\bibitem{r7}
J.S.~{Purcell}, \emph{{Hyperbolic Knot Theory}}, {\emph{arXiv e-prints} (2020)
  arXiv:2002.12652} [\href{https://arxiv.org/abs/2002.12652}{{\ttfamily
  2002.12652}}].

\bibitem{r8}
S.~{Gukov}, \emph{{Three-Dimensional Quantum Gravity, Chern-Simons Theory, and
  the A-Polynomial}},
  \href{https://doi.org/10.1007/s00220-005-1312-y}{\emph{Communications in
  Mathematical Physics} {\bfseries 255} (2005) 577}
  [\href{https://arxiv.org/abs/hep-th/0306165}{{\ttfamily hep-th/0306165}}].

\bibitem{r10}
E.~Witten, \emph{Quantization of {Chern-Simons} gauge theory with complex gauge
  group}, \href{https://doi.org/10.1007/BF02099116}{\emph{Communications in
  Mathematical Physics} {\bfseries 137} (1991) 29}.

\bibitem{r11}
E.~Witten, \emph{Quantum field theory and the {Jones} polynomial},
  \href{https://doi.org/10.1007/bf01217730}{\emph{Communications in
  Mathematical Physics} {\bfseries 121} (1989) 351}.

\bibitem{r13}
R.K.~{Kaul}, \emph{{Chern-Simons theory, coloured-oriented braids and link
  invariants}}, \href{https://doi.org/10.1007/BF02102019}{\emph{Communications
  in Mathematical Physics} {\bfseries 162} (1994) 289}
  [\href{https://arxiv.org/abs/hep-th/9305032}{{\ttfamily hep-th/9305032}}].

\bibitem{r12}
V.~{Balasubramanian}, J.R.~{Fliss}, R.G.~{Leigh} and O.~{Parrikar},
  \emph{{Multi-boundary entanglement in Chern-Simons theory and link
  invariants}}, \href{https://doi.org/10.1007/JHEP04(2017)061}{\emph{Journal of
  High Energy Physics} {\bfseries 2017} (2017) 61}
  [\href{https://arxiv.org/abs/1611.05460}{{\ttfamily 1611.05460}}].

\bibitem{r20}
J.R.~{Fliss}, \emph{{Knots, links, and long-range magic}},
  \href{https://doi.org/10.1007/JHEP04(2021)090}{\emph{Journal of High Energy
  Physics} {\bfseries 2021} (2021) 90}
  [\href{https://arxiv.org/abs/2011.01962}{{\ttfamily 2011.01962}}].

\bibitem{r14}
R.~Kaul and T.~Govindarajan, \emph{Three-dimensional {Chern-Simons} theory as a
  theory of knots and links: (ii). multicoloured links},
  \href{https://doi.org/https://doi.org/10.1016/0550-3213(93)90251-J}{\emph{Nuclear
  Physics B} {\bfseries 393} (1993) 392}.

\bibitem{r15}
R.K.~{Kaul}, \emph{{Chern-Simons Theory, Knot Invariants, Vertex Models and
  Three-manifold Invariants}},
  [\href{https://arxiv.org/abs/hep-th/9804122}{{\ttfamily hep-th/9804122}}].

\bibitem{r16}
R.K.~{Kaul}, \emph{{Complete Solution of SU(2) Chern-Simons Theory}},
  {\emph{arXiv e-prints} (1992) hep}
  [\href{https://arxiv.org/abs/hep-th/9212129}{{\ttfamily hep-th/9212129}}].

\bibitem{r17}
V.~{Balasubramanian}, M.~{DeCross}, J.~{Fliss}, A.~{Kar}, R.G.~{Leigh} and
  O.~{Parrikar}, \emph{{Entanglement entropy and the colored Jones
  polynomial}}, \href{https://doi.org/10.1007/JHEP05(2018)038}{\emph{Journal of
  High Energy Physics} {\bfseries 2018} (2018) 38}
  [\href{https://arxiv.org/abs/1801.01131}{{\ttfamily 1801.01131}}].

\bibitem{r35}
S.~{Bai}, \emph{{Hyperbolic Brunnian links}}, {\emph{arXiv e-prints} (2021)
  arXiv:2104.12637} [\href{https://arxiv.org/abs/2104.12637}{{\ttfamily
  2104.12637}}].

\bibitem{r18}
S.~{Dwivedi}, V.K.~{Singh}, S.~{Dhara}, P.~{Ramadevi}, Y.~{Zhou} and
  L.K.~{Joshi}, \emph{{Entanglement on linked boundaries in Chern-Simons theory
  with generic gauge groups}},
  \href{https://doi.org/10.1007/JHEP02(2018)163}{\emph{Journal of High Energy
  Physics} {\bfseries 2018} (2018) 163}
  [\href{https://arxiv.org/abs/1711.06474}{{\ttfamily 1711.06474}}].

\bibitem{r19}
S.~{Dwivedi}, V.K.~{Singh}, P.~{Ramadevi}, Y.~{Zhou} and S.~{Dhara},
  \emph{{Entanglement on multiple S $^{2}$ boundaries in Chern-Simons theory}},
  \href{https://doi.org/10.1007/JHEP08(2019)034}{\emph{Journal of High Energy
  Physics} {\bfseries 2019} (2019) 34}
  [\href{https://arxiv.org/abs/1906.11489}{{\ttfamily 1906.11489}}].

\bibitem{r24}
H.~{Ooguri} and C.~{Vafa}, \emph{{Knot invariants and topological strings}},
  \href{https://doi.org/10.1016/S0550-3213(00)00118-8}{\emph{Nuclear Physics B}
  {\bfseries 577} (2000) 419}
  [\href{https://arxiv.org/abs/hep-th/9912123}{{\ttfamily hep-th/9912123}}].

\bibitem{r25}
J.M.F.~{Labastida}, M.~{Mari{\~n}o} and C.~{Vafa}, \emph{{Knots, links and
  branes at large N}},
  \href{https://doi.org/10.1088/1126-6708/2000/11/007}{\emph{Journal of High
  Energy Physics} {\bfseries 2000} (2000) 007}
  [\href{https://arxiv.org/abs/hep-th/0010102}{{\ttfamily hep-th/0010102}}].

\bibitem{r26}
J.M.F.~{Labastida} and M.~{Marino}, \emph{{A New Point of View in the Theory of
  Knot and Link Invariants}}, {\emph{arXiv Mathematics e-prints} (2001)
  math/0104180} [\href{https://arxiv.org/abs/math/0104180}{{\ttfamily
  math/0104180}}].

\bibitem{r31}
J.M.~{Isidro}, J.M.F.~{Labastida} and A.V.~{Ramallo}, \emph{{Coset
  constructions in Chern-Simons gauge theory}},
  \href{https://doi.org/10.1016/0370-2693(92)90480-R}{\emph{Physics Letters B}
  {\bfseries 282} (1992) 63}
  [\href{https://arxiv.org/abs/hep-th/9201027}{{\ttfamily hep-th/9201027}}].

\bibitem{r32}
J.M.~{Isidro}, J.M.F.~{Labastida} and A.V.~{Ramallo}, \emph{{Polynomials for
  torus links from Chern-Simons gauge theories}},
  \href{https://doi.org/10.1016/0550-3213(93)90632-Y}{\emph{Nuclear Physics B}
  {\bfseries 398} (1993) 187}
  [\href{https://arxiv.org/abs/hep-th/9210124}{{\ttfamily hep-th/9210124}}].

\bibitem{r36}
D.~Friedan, Z.~Qiu and S.~Shenker, \emph{Conformal invariance, unitarity, and
  critical exponents in two dimensions},
  \href{https://doi.org/10.1103/PhysRevLett.52.1575}{\emph{Phys. Rev. Lett.}
  {\bfseries 52} (1984) 1575}.

\bibitem{r30a}
P.~{Ramadevi}, T.R.~{Govindarajan} and R.K.~{Kaul}, \emph{{Knot invariants from
  rational conformal field theories}},
  \href{https://doi.org/10.1016/0550-3213(94)00102-2}{\emph{Nuclear Physics B}
  {\bfseries 422} (1994) 291}
  [\href{https://arxiv.org/abs/hep-th/9312215}{{\ttfamily hep-th/9312215}}].

\bibitem{r34}
A.~{Kitaev}, \emph{{Anyons in an exactly solved model and beyond}},
  \href{https://doi.org/10.1016/j.aop.2005.10.005}{\emph{Annals of Physics}
  {\bfseries 321} (2006) 2}
  [\href{https://arxiv.org/abs/cond-mat/0506438}{{\ttfamily
  cond-mat/0506438}}].

\bibitem{r22}
Z.P.~{Ramadevi}, \emph{{SU(N) quantum Racah coefficients \& non-torus links}},
  {\emph{arXiv e-prints} (2011) arXiv:1107.3918}
  [\href{https://arxiv.org/abs/1107.3918}{{\ttfamily 1107.3918}}].

\bibitem{r21}
S.G.~{Naculich} and H.J.~{Schnitzer}, \emph{{Level-rank duality of the U(N) WZW
  model, Chern-Simons theory, and 2d qYM theory}},
  \href{https://doi.org/10.1088/1126-6708/2007/06/023}{\emph{Journal of High
  Energy Physics} {\bfseries 2007} (2007) 023}
  [\href{https://arxiv.org/abs/hep-th/0703089}{{\ttfamily hep-th/0703089}}].

\bibitem{r23}
T.~{Dimofte} and S.~{Gukov}, \emph{{Quantum Field Theory and the Volume
  Conjecture}}, {\emph{arXiv e-prints} (2010) arXiv:1003.4808}
  [\href{https://arxiv.org/abs/1003.4808}{{\ttfamily 1003.4808}}].

\bibitem{r27}
S.~{Zhu}, \emph{{New structures for colored HOMFLY-PT invariants}},
  {\emph{arXiv e-prints} (2021) arXiv:2105.02037}
  [\href{https://arxiv.org/abs/2105.02037}{{\ttfamily 2105.02037}}].

\bibitem{r28}
Q.~{Chen}, K.~{Liu}, P.~{Peng} and S.~{Zhu}, \emph{{Congruent skein relations
  for colored HOMFLY-PT invariants and colored Jones polynomials}},
  {\emph{arXiv e-prints} (2014) arXiv:1402.3571}
  [\href{https://arxiv.org/abs/1402.3571}{{\ttfamily 1402.3571}}].

\bibitem{r29}
S.~{Gukov} and M.~{Stosic}, \emph{{Homological algebra of knots and BPS
  states}}, {\emph{arXiv e-prints} (2011) arXiv:1112.0030}
  [\href{https://arxiv.org/abs/1112.0030}{{\ttfamily 1112.0030}}].

\bibitem{r30}
E.~{Gorsky}, S.~{Gukov} and M.~{Stosic}, \emph{{Quadruply-graded colored
  homology of knots}}, {\emph{arXiv e-prints} (2013) arXiv:1304.3481}
  [\href{https://arxiv.org/abs/1304.3481}{{\ttfamily 1304.3481}}].

\bibitem{r33}
N.~{Kubo} and S.~{Yokoyama}, \emph{{Topological phase, spin Chern-Simons theory
  and level rank duality on lens space}}, {\emph{arXiv e-prints} (2021)
  arXiv:2108.09300} [\href{https://arxiv.org/abs/2108.09300}{{\ttfamily
  2108.09300}}].

\end{thebibliography}\endgroup
\end{document}